\documentclass[12pt,letter,a4paper]{article}
\usepackage{amsmath}
\usepackage{amsfonts}
\begin{document}
\title{Extending the Calculus of Moving Surfaces to Higher Orders}
\author{Keith C. Afas\thanks{Corresponding Author}  \thanks{kafas@uwo.ca}\\ \\ University of Western Ontario\\ London, Ontario, N6A 5B7, CANADA}
\maketitle
\tableofcontents
\section*{Abstract}
\textrm{In 2010, a book published on the work of Jaques Hadamard, entitled "Introduction to Tensor Analysis and the Calculus of Moving Surfaces'' by Dr. Pavel Grinfeld, proposed an extension of Hadamard's work to ultimately allow principles of tensorial invariance (introduced in the now ubiquitous Tensor Calculus) to be raised and extended to notions of time-dependent and moving surfaces \cite{CMS}. Coined the "Calculus of Moving Surfaces'' (CMS), notions of Invariant Time Derivatives on Arbitrary Surface/Ambient Tensors, Surface Velocities, and time derivatives of time-dependent Volume \& Surface Integrals were introduced. This discipline was shown to have diverse application in problems such as shape optimization, boundary value problems, and time evolutions of surfaces. This paper focuses on extending Invariant Time Derivatives to other Surface Objects and to Higher Orders \cite{CMS}.}
\section{Introduction}
\textrm{When Considering Manifolds in General, it is established that the field of Differential Geometry governs the Spatial Connections defined on a Manifold \cite{DifGeo1}, }$\Sigma$\textrm{ and establishes rules for the Definition and Manipulation of Tensors defined on sections of the Tangent Bundle of }$\Sigma$\textrm{ denoted by }$T\Sigma$\textrm{ \cite{TenAn1, DifGeo1}. Tensors on the manifold usually map regions of the Arithmetic Space of }$\Sigma$\textrm{ which is just a region of }$\mathbb{E}^2$\textrm{ denoted by }$S$\textrm{ to the space of all tensors of rank }$(n,m)$\textrm{ defined on }$\Sigma$\textrm{ notated by }$\tilde{\mathcal{T}}^{(n,m)}(\Sigma)$\textrm{ \cite{DifGeo1}. This is encapsulated by:}
\begin{equation}
F:S\to\tilde{\mathcal{T}}^{(n,m)}(\Sigma)
\end{equation}
\textrm{where }$\tilde{\mathcal{T}}^{(n,m)}(\Sigma)$\textrm{ is usually explained as a space formed by Cartesian Products and Exterior Products on }$T\Sigma$\textrm{.}
\\\\\phantom{idnt}
\textrm{When the manifold is parametrized by a parameter }$t$\textrm{ (commonly identified by "time"), the Manifold is commonly denoted by }$\Sigma_t$\textrm{; the parameter }$t$\textrm{ usually varies in a smooth way such that there are no 'jumps' between }$\Sigma_t$\textrm{ and }$\Sigma_{t+\epsilon}$\textrm{ (where }$\epsilon$\textrm{ is an infinitesimal). All Tensors defined on }$\Sigma_t$\textrm{ usually map the Cartesian product of the Arithmetic Space with a closed and compact interval of the Real Numbers representing the range of time in the following manner:}
\begin{equation}
F:S\times\mathbb{R}\to\tilde{\mathcal{T}}^{(n,m)}(\Sigma_t)
\end{equation}
\textrm{ In general, Tensors are usually constructed as being scalar multiples of the Basis of the Manifold they reside on and are always decomposable into the Basis, and Components with respect to this Basis \cite{TenAn2}. When the Coordinate System undergoes a Change of Basis, in order to preserve the Geometrical Meaning of a Tensor, its Basis undergoes a compensatory transformation, as does its components in a manner which is dependent on the Jacobian of the Change of Basis.}\\\\\phantom{idnt}
\textrm{When time is introduced, an ambiguity is also introduced with Change of Bases. Grinfeld\cite{CMS, BetterCMS} established that for a }$\textbf{Dynamic}$\textrm{ 2-Manifold }$\Sigma_t$\textrm{, when the Coordinate System undergoes a passive transformation, operators which usually preserved the Jacobian-dependent transformation of their tensorial arguments on }$\Sigma$\textrm{ no longer produced tensors }$\Sigma_t$\textrm{ such as the partial time derivative }$\partial_t$\textrm{ \cite{CMS}.}\\\\\phantom{idnt}
\textrm{In order to recover objects which would transform in identical jacobian dependent ways on static and dynamic surfaces, Grinfeld developed the field of the Calculus of Moving Surfaces (CMS). This paper will explore the implications brought up by Grinfeld's formulation of CMS and introduce more objects which are useful metrics for the Motion of Dynamic 2-Manifolds as well as some tensorial relations which are present on Dynamic 2-Manifolds.}\\\\\phantom{idnt}
\textrm{Two Jewels of CMS are the Surface Velocity, }$\tilde{C}$\textrm{ which is a scalar tensorial field that serves as a measurement of how the 2-Manifold infinitesimally moves in the normal direction, and the definition of the Tensorial Time Derivative Operator, }$\dot{\nabla}$\textrm{ defined for Tensors }$F\in\tilde{\mathcal{T}}^{(n,m)}(\Sigma_t)$\textrm{\cite{CMS}.}
\subsection{The Normal Surface Velocity}
\textrm{The Surface Velocity is a scalar field parameter of a surface which describes the normal speed of the Dynamic Surface, defined as:}
\begin{equation}
\tilde{C}=\textbf{V}\cdot\textbf{N}
\end{equation}
\textrm{Where }$\textbf{V}$\textrm{ is the Time Partial Derivative of the Embedding of the 2-Manifold and }$\textbf{N}$\textrm{ is the Unit Normal Field defined on the 2-Manifold. It is unique in that its contributory tensor }$\textbf{V}$\textrm{ is not a Tensor, but it is a parameter which is invariant with respect to changes to the Surface/Ambient set of coordinates \cite{CMS}.}
\subsection{Introduction to the Invariant Time Derivative}
\textrm{The Invariant Time Derivative is similar to the Surface Velocity in that it too, is an operator which remains invariant under Surface/Ambient Coordinate Changes on Invariant Fields of rank (0,0). As such, it is an appropiate choice of a temporal derivative on a surface.}\\\\\phantom{idnt}
\textrm{In general, this operator is unique in that it allows for the preservation of the scalar invariant nature of its operand, preserves its rank, and captures the rate that the scalar field changes on a surface \cite{CMS}. It is defined to satisfied four axioms:}
\begin{description}
\item[I.]{The Principle of Invariance: }$\dot{\nabla}:\tilde{\mathcal{T}}^{(0,0)}(\Sigma_t)\to\tilde{\mathcal{T}}^{(0,0)}(\Sigma_t)$
\item[II.]{Associativity with Addition: }$\dot{\nabla}(\psi+\sigma)=\dot{\nabla}\psi+\dot{\nabla}\sigma$
\item[III.]{Liebniz's Rule: }$\dot{\nabla}(\psi\sigma)=\psi\dot{\nabla}\sigma+\sigma\dot{\nabla}\psi$
\item[IV.]{Spatiotemporal Semi-Commutation: }$(\dot{\nabla}\nabla_\alpha-\nabla_\alpha\dot{\nabla}-\tilde{C}B^\beta_\alpha\nabla_\beta)\psi=0$
\end{description}
\textrm{With the 1st Axiom being the definitive one. Grinfeld generalizes the concept of an Invariant Time Derivative operating on an Invariant field to the }$\textbf{Tensorial Time Derivative}$\textrm{ which also behaves as a Tensor of a rank identical to the tensor it is applied to (ie. the Tensorial Time Derivative of a tensor of rank (n,m) will transform under coordinate changes as a tensor of rank (n,m) would) \cite{CMS}. This Time Operator reduces to the Invariant Time Derivative since for a tensor of rank (0,0), the time derivative will transform as a tensor of rank (0,0) would, which is not at all (ie. ``invariant''). This new operator is defined to satisfy six axioms:}
\begin{description}
\item[I.]{The Principle of Tensoriality: }$\dot{\nabla}:\tilde{\mathcal{T}}^{(n,m)}(\Sigma_t)\to\tilde{\mathcal{T}}^{(n,m)}(\Sigma_t)$
\item[II.]{Associativity with Addition: }$\dot{\nabla}(\psi+\sigma)=\dot{\nabla}\psi+\dot{\nabla}\sigma$
\item[III.]{Liebniz's Rule: }$\dot{\nabla}(\psi\otimes\sigma)=\psi\otimes\dot{\nabla}\sigma+\dot{\nabla}\psi\otimes\sigma$
\item[IV.]{Metrilinicity with Respect to the Basis of the Ambient Space: }$\dot{\nabla}\textbf{Z}_i=\textbf{0}$
\item[VI.]{Geometric Fixing Principle: }$\dot{\nabla}\textbf{R}=\tilde{C}\textbf{N}$
\item[VII.]{Invariant Semi-Commutation:}$(\dot{\nabla}\nabla_\alpha-\nabla_\alpha\dot{\nabla}-\tilde{C}B^\beta_\alpha\nabla_\beta)\psi=0$
\end{description}
\textrm{Again, where the 1st Axiom is the definitive one. Informally the Invariant Time Derivative, and more generally, the Tensorial Time Derivative, captures the temporal rate of change of a general tensor field defined on a dynamic surface. For a static surface, the operator behaves like a time partial derivative. Thus it satisfies the same analogue as the Covariant Derivative, defined on the connection of }$\Sigma$\textrm{ or the Ambient Space it is embedded into \cite{DifGeo1}.}
\subsection{The Invariant Time Derivative}
\textrm{The expression of the Invariant Time Derivative is shown in Local Coordinates. Considering a general tensor }$\psi\in\tilde{\mathcal{T}}^{(0,0)}(\Sigma_t)$\textrm{ such that }$\psi=\psi(S,t)$\textrm{ defined on an arbitrary dynamic 2-Manifold embedded by the following ambient vector }$\textbf{R}(S,t)$\textrm{ where the surface coordinates, }$S=(S^1,S^2)$\textrm{, are indexed as }$S^\alpha$\textrm{, by definition the Invariant Time Derivative is defined as:}
\begin{equation} \label{DerDef}
(\dot{\nabla})\psi=\left(\partial_t-V^\alpha\nabla_\alpha\right)\psi
\end{equation}
\textrm{Where }$\nabla_\alpha$\textrm{ indicates the surface covariant derivative defined on the Connection shared between }$\Sigma$\textrm{ \& }$\Sigma_t$\textrm{, }$\partial_t$\textrm{ indicates partial differentiation with respect to time, and }$V^\alpha$\textrm{ is defined as }$V^\alpha=\textbf{S}^\alpha\cdot\textbf{V}$\textrm{, where }$\textbf{S}^\alpha$\textrm{ is the dual basis to }$\Sigma_t$\textrm{, and }$\textbf{V}$\textrm{ is defined as:}
\begin{equation*}
\textbf{V}\in\mathbb{V}^3\phantom{.},\phantom{.}\textbf{V}=\partial_t\textbf{R}
\end{equation*}
\textrm{Where }$\mathbb{V}^3$\textrm{ is defined as the set of Vectors formable from the Euclidean 3-space, }$\mathbb{R}^3$\textrm{. This is extended to tensors in }$\tilde{\mathcal{T}}^{(n,m)}(\Sigma_t)$\textrm{ by considering an Invariant Tensor constructed from its tensorial components contracted to the basis of }$\Sigma_t$\textrm{. The Full Expression can be substituted into The Definition of the Invariant Time Derivative for Invariant Tensors and expanded \& collected accordingly. Considering a Surface Invariant Tensor, }$\textbf{T}\in\tilde{\mathcal{T}}^{(1,0)}(\Sigma_t)$\textrm{ defined using }$\Sigma_t$\textrm{'s basis as }$\textbf{T}=T^\alpha\textbf{S}_\alpha$\textrm{; its Invariant Time Derivative is expressed as:}
\begin{equation*}
\dot{\nabla}\textbf{T}=\partial_t\textbf{T}-V^\beta\nabla_\beta\textbf{T}\phantom{..}\to\phantom{..}\dot{\nabla}\left(T^\alpha\textbf{S}_\alpha\right)=\partial_t\left(T^\alpha\textbf{S}_\alpha\right)-V^\beta\nabla_\beta\left(T^\alpha\textbf{S}_\alpha\right)
\end{equation*}
\textrm{Since ordinary time partial derivatives are well defined for tensors, as are surface covariant derivatives; and using the obtainable fact that the invariant time derivative obeys the product rule, then the time derivatives can be extended to tensors of arbitrary classification and rank by expanding the left side and the right side of the above equation, and grouping/matching terms.}
\subsection{The Tensorial Time Derivative}
\textrm{Consider a tensor field defined on a dynamic surface, }$\tilde{\phi}\in\tilde{\mathcal{T}}^{(2,2)}(\Sigma_t)$\textrm{ containing one contra- and one co-variant ambient index, and containing one contra- and one co-variant surface index. The form of the invariant surface field is given as:}
\begin{equation*}
\tilde{\phi}=\phi^{i\alpha}_{j\beta}(\textbf{Z}_i\otimes\textbf{Z}^j\otimes\textbf{S}_\alpha\otimes\textbf{S}^\beta)
\end{equation*}
\textrm{The form of the invariant time derivative of such a multi-indicial tensor's components is given by:}
\begin{equation}
\dot{\nabla}\phi^{i\alpha}_{j\beta}=\left(\partial_t-V^\gamma\nabla_\gamma\right)\phi^{i\alpha}_{j\beta}+V^k\left(\Gamma^i_{mk}\phi^{m\alpha}_{j\beta}-\Gamma^n_{jk}\phi^{i\alpha}_{n\beta}\right)+\dot{\Gamma}^\alpha_\sigma\phi^{i\sigma}_{j\beta}-\dot{\Gamma}^\lambda_\beta\phi^{i\alpha}_{j\lambda}
\end{equation}
\textrm{Where the Christoffel Time Symbol, }$\dot{\Gamma}^\alpha_\beta$\textrm{ is defined as:}
\begin{equation}\label{ChrisTime}
\dot{\Gamma}^\alpha_\beta=\nabla_\beta V^\alpha-\tilde{C}B^\alpha_\beta
\end{equation}
\textrm{Taking }$B^\alpha_\beta$\textrm{ to be the Curvature Tensor defined as:}
\begin{equation}
B^\alpha_\beta=-\textbf{S}^\alpha\cdot\nabla_\beta\textbf{N}=\textbf{N}\cdot\nabla_\beta\textbf{S}^\alpha
\end{equation}
\textrm{Taking these, the results on various differential objects residing on dynamic surfaces can be discussed further. Thus, the Tensorial Time Derivative operates on various Surface objects as per the following \cite{CMS}:}
\begin{equation}\label{CMStable}
\left.\begin{matrix}
\dot{\nabla}\textbf{Z}^i=\dot{\nabla}\textbf{Z}_i=\textbf{0} \\
\dot{\nabla}Z^{ij}=\dot{\nabla}Z_{ij}=0 \\
\dot{\nabla}\textbf{S}_\alpha=\textbf{N}\nabla_\alpha\tilde{C}\\
\dot{\nabla}\textbf{S}^\alpha=\textbf{N}\nabla^\alpha\tilde{C}\\
\dot{\nabla}S^{\alpha\beta}=\dot{\nabla}S_{\alpha\beta}=0\\
\dot{\nabla}\textbf{N}=-\textbf{S}_\alpha\nabla^\alpha\tilde{C}=-\textbf{S}^\alpha\nabla_\alpha\tilde{C} \\
\dot{\nabla}B^{\alpha}_{\beta}=\nabla_\beta\nabla^\alpha\tilde{C}+B^\alpha_\gamma B^\gamma_\beta\tilde{C} \\
\dot{\nabla}\epsilon^{ijk}=\dot{\nabla}\epsilon_{ijk}=0 \\
\dot{\nabla}\epsilon^{\alpha\beta}=\dot{\nabla}\epsilon_{\alpha\beta}=0\\
\dot{\nabla}\delta^i_{\phantom{.}j}=\dot{\nabla}\delta^\alpha_{\phantom{.}\beta}=0
\end{matrix}\right.
\end{equation}
\textrm{Where }$\epsilon^{ijk}$\textrm{ \& }$\epsilon_{ijk}$\textrm{, and }$\epsilon^{\alpha\beta}$\textrm{ and }$\epsilon_{\alpha\beta}$\textrm{ are the Ambient \& Surface Levi-Civita symbols, respectively. Also, }$\delta^i_{\phantom{.}j}$\textrm{ and }$\delta^\alpha_{\phantom{.}\beta}$\textrm{ are the Ambient and Surface Kronecker-Delta symbols, respectively. Though this is a fairly comprehensive list, there are several critical relations which are not exhibited, such as }$\dot{\nabla}C$\textrm{, as well as no expression extending the above equations to second orders derivatives (ie. }$\dot{\nabla}^2\textbf{N}$\textrm{). Just as second partial derivatives have found use in Physics so should second order Tensorial Time Derivatives find use, if not more use than the standard second order partial time derivative }$\partial_t^2$\textrm{.}
\section{Second-Order Tensorial Time-Derivative}
\subsection{Introduction to the Second Order Invariant Time-Derivative}
\textrm{The Second-Order Tensorial Time-Derivative }$\dot{\nabla}^2$\textrm{ can be applied to tensors of rank (0,0). Considering the tensor }$\psi\in\tilde{\mathcal{T}}^{(0,0)}(\Sigma_t)$\textrm{ defined on }$\Sigma_t$\textrm{, }$\dot{\nabla}\psi$\textrm{ must also be an element of }$\tilde{\mathcal{T}}^{(0,0)}(\Sigma_t)$\textrm{; also iterating the expression again, }$\dot{\nabla}(\dot{\nabla}\psi)$\textrm{ must be a an element of }$\tilde{\mathcal{T}}^{(0,0)}(\Sigma_t)$\textrm{. So the Second Order Invariant Time Derivative can be defined as the Invariant Time Derivative Operator iterated twice: }$\dot{\nabla}^2\psi=\dot{\nabla}(\dot{\nabla}\psi)$\textrm{; it is a unique operator which preserves tensors and satisfies the following property:}
\begin{equation}
\dot{\nabla}^2:\tilde{\mathcal{T}}^{(0,0)}(\Sigma_t)\to\tilde{\mathcal{T}}^{(0,0)}(\Sigma_t)
\end{equation}
\textrm{In general the operator can be expressed in local coordinates, but the function of the operator essentially }$\textbf{preserves}$\textrm{ }$\tilde{\mathcal{T}}^{(0,0)}(\Sigma_t)$\textrm{ tensors defined on }$\Sigma_t$\textrm{. This concept can be extended to a Second Order }\textbf{Tensorial}\textrm{ Time Derivative Operator which serves this function for the Second-Order iteration of the Operator applied to tensors of rank greater than }$(0,0)$\textrm{:}
\begin{equation}
\dot{\nabla}^2:\tilde{\mathcal{T}}^{(n,m)}(\Sigma_t)\to\tilde{\mathcal{T}}^{(n,m)}(\Sigma_t)
\end{equation}
\subsection{Local Expression of the Operator}
\textrm{The expression of the Second Order Invariant Time Derivative Operator acting on rank (0,0) Tensors can be expressed in local coordinates. If the operator is iterated twice, the twice-iterated operator takes the form of:}
\begin{equation*}
\dot{\nabla}^2\psi=\dot{\nabla}(\dot{\nabla}\psi)
\end{equation*}
\textrm{Applying the operator from outside inwards \& and expanding,}
\begin{equation*}
\dot{\nabla}^2\psi=\partial_t(\dot{\nabla}\psi)-V^\alpha\partial_\alpha(\dot{\nabla}\psi)
\end{equation*}
\textrm{Applying the operator once more inwards, the form of the twice-iterated operator is as follows:}
\begin{equation}
\dot{\nabla}^2\psi=\partial_t\left(\partial_t\psi-V^\alpha\partial_\alpha\psi\right)-V^\alpha\partial_\alpha\left(\partial_t\psi-V^\beta\partial_\beta\psi\right)
\end{equation}
\textrm{After expanding and recognizing the symmetry of the partial derivatives, }$\partial_t\partial_\alpha-\partial_\alpha\partial_t=0$\textrm{, the operator takes the following form:}
\begin{equation}\label{SOTTD}
\dot{\nabla}^2\psi=\partial_t^2\psi-\left(\partial_tV^\alpha-V^\beta \partial_\beta V^\alpha\right)\partial_\alpha\psi-V^\alpha\left(2\partial_t\partial_\alpha\psi-V^\beta\partial_\beta\partial_\alpha\psi\right)
\end{equation}
\textrm{This is the Second Invariant Time Derivative Operator \& can be extended to Tensors of Higher Order.}
\subsection{Applications of the Second Order Tensorial Time Derivative}
\textrm{If one knows the Second-Order Invariant Time-Derivative of }$\psi$\textrm{, denoted by }$\dot{\nabla}^2\psi$\textrm{, then one can obtain its Second-Order time partial derivative }$\partial_t^2$\textrm{ and vice-versa. The local coordinate expression can also be useful for using the Second-Order time partial derivative }$\partial_t^2$\textrm{ of commonly used differential objects (}$\textbf{N},\textbf{S}_\alpha,etc.$\textrm{) since their SOTTD can easily be found by the table in Eq.(\ref{CMStable}).}\\\\\phantom{idnt}	
\textrm{It is worth noting that the SOTTD for arbitrary tensors of higher rank must be re-calculated and are }\textbf{not}\textrm{ contained in Eq.(\ref{SOTTD})}
\section{The Temporal Curvature Tensor}
\subsection{Motivation for the Temporal Curvature Tensor}
\textrm{Analyzing }$\Sigma$\textrm{'s connections, the Covariant Derivatives }$\nabla_\alpha$\textrm{ }\textbf{do not}\textrm{ commute \cite{CMS,TenAn2,DifGeo1} on Tensors of rank greater than (0,0). This lack of commutation is given by the Riemann Curvature Tensor }$R^\delta_{\gamma\alpha\beta}$\textrm{. However, the Invariant Time Derivative }$\dot{\nabla}$\textrm{ contain connection properties of its own \cite{CMS}. Since it semi-commutes on an invariant, the question naturally arises on whether it has a special commutation tensor, }$\dot{R}^\beta_\alpha\gamma$\textrm{ which appears after commuting the Tensorial Time Derivative }$\dot{\nabla}$\textrm{, and the Covariant Derivative }$\nabla_\alpha$\textrm{. On tensors of rank greater than (0,0)}\\\\\phantom{idnt}
\textrm{Often when analyzing two linear differential operators, }$O$\textrm{ \& }$\Omega$\textrm{, it is often inquired to the effect of their commutation on a tensor field \cite{CommOps}, }$\mathcal{T}$\textrm{ by defining a new operator, their commutator:}
\begin{equation}
(\kappa)\mathcal{T}=\left(O\Omega-\Omega O\right)\mathcal{T}
\end{equation}
\textrm{The commutation relation can obviously depend on O \& }$\Omega$\textrm{, but in some cases, the commutation can take different forms depending on the rank of }$\mathcal{T}$\textrm{. For example, multiplication is an operator well defined for Scalars and is known to commute for scalars (ie. rank (0,0) tensors). But in the case that one operator is the Cross Product, and the other operator is the standard Dot Product, then result will greatly depend on the order of the Tensor Field and, more specifically, its components.}\\\\\phantom{idnt}
\textrm{A famous example is the commutation operator }$(\nabla_\alpha\nabla_\beta-\nabla_\beta\nabla_\alpha)$\textrm{. When applied to a tensor of rank (0,0), the operators commute \cite{CMS}, and so the relation exists as:}
\begin{equation}\label{SurfComm}
\left(\nabla_\alpha\nabla_\beta-\nabla_\beta\nabla_\alpha\right)\psi=0
\end{equation}
\textrm{However, when the form of the tensor field is chosen to have one contravariant index, such as }$\mathcal{T}=\psi^\alpha\textbf{S}_\alpha$\textrm{, then the commutation relation shows as \cite{TenAn2, DifGeo1}:}
\begin{equation}\label{SurfCommRiem}
\left(\nabla_\alpha\nabla_\beta-\nabla_\beta\nabla_\alpha\right)\psi^\gamma=R^\gamma_{\phantom{.}\delta\alpha\beta}\psi^\delta\phantom{.},\phantom{.}\left(\nabla_\alpha\nabla_\beta-\nabla_\beta\nabla_\alpha\right)\psi_\gamma=-R^\delta_{\phantom{.}\gamma\alpha\beta}\psi_\delta
\end{equation}
\textrm{Where }$R^\gamma_{\phantom{.}\delta\alpha\beta}$\textrm{ is the Riemann Curvature Tensor for the Surface. Using the commutation operator }$(\dot{\nabla}\nabla_\alpha-\nabla_\alpha\dot{\nabla})$\textrm{ on scalar fields defined on manifolds results in the following relation \cite{CMS}:}
\begin{equation*}
(\dot{\nabla}\nabla_\alpha-\nabla_\alpha\dot{\nabla})\psi=\tilde{C}B^\beta_\alpha\nabla_\beta\psi
\end{equation*}
\textrm{Equally valid in the form as:}
\begin{equation}
(\dot{\nabla}\nabla_\alpha-\nabla_\alpha\dot{\nabla}-\tilde{C}B^\beta_\alpha\nabla_\beta)\psi=0
\end{equation}
\textrm{Or the commutation form as:}
\begin{equation}\label{SpaTemInv}
([\dot{\nabla},\nabla_\alpha]-\tilde{C}B^\beta_\alpha\nabla_\beta)\psi=0
\end{equation}
\textrm{This equation describes  very similar statement as Eq.(\ref{SurfComm}) Since by Eq.(\ref{DerDef}), it is established that the Invariant Time Derivative Operator contains a covariant derivative, the above relations beg the question: Are there any relations pertaining to the following commutation operator?}
\begin{equation*}
(\dot{\nabla}\nabla_\alpha-\nabla_\alpha\dot{\nabla})\psi^\beta
\end{equation*}
\subsection{Derivation of the Temporal Curvature Tensor}
\textrm{The expansion of the operator begins, obeying the tensorial conventions:}
\begin{equation*}
(\dot{\nabla}\nabla_\alpha-\nabla_\alpha\dot{\nabla})\psi^\beta=\dot{\nabla}(\nabla_\alpha\psi^\beta)-\nabla_\alpha(\dot{\nabla}\psi^\beta)
\end{equation*}
\textrm{The right side is expanded from the outside inwards recalling the use of the Time Christoffel Symbols, }$\dot{\Gamma}^\alpha_\beta$\textrm{ when the Tensorial Time Derivative is applied to tensors,}
\begin{equation*}
\partial_t(\nabla_\alpha\psi^\beta)-V^\gamma\nabla_\gamma\nabla_\alpha\psi^\beta+\dot{\Gamma}^\beta_\gamma\nabla_\alpha\psi^\gamma-\dot{\Gamma}^\gamma_\alpha\nabla_\gamma\psi^\beta-\partial_\alpha(\dot{\nabla}\psi^\beta)-\Gamma^\beta_{\alpha\gamma}\dot{\nabla}\psi^\gamma \\
\end{equation*}
\textrm{Expanding the first }$\nabla_\alpha\psi^\beta$\textrm{ and the first }$\dot{\nabla}\psi^\beta$\textrm{ further inwards and recalling that the partial derivatives, }$\partial_t$\textrm{ \& }$\partial_\alpha$\textrm{ commute, the right side can be further simplified to the following:}
\begin{equation*}
\partial_t(\Gamma^\beta_{\alpha\gamma}\psi^\gamma)-V^\gamma\nabla_\gamma\nabla_\alpha\psi^\beta+\dot{\Gamma}^\beta_\gamma\nabla_\alpha\psi^\gamma-\dot{\Gamma}^\gamma_\alpha\nabla_\gamma\psi^\beta+\partial_\alpha(V^\gamma\nabla_\gamma\psi^\beta-\dot{\Gamma}^\beta_\gamma\psi^\gamma)-\Gamma^\beta_{\alpha\gamma}\dot{\nabla}\psi^\gamma \\
\end{equation*}
\textrm{This expression after expanding and contracting identities can be re-arranged to the following:}
\begin{equation*}
\psi^\gamma\partial_t\Gamma^\beta_{\alpha\gamma}-\Gamma^\beta_{\alpha\gamma}(\dot{\nabla}\psi^\gamma-\partial_t\psi^\gamma)-V^\gamma(\nabla_\gamma\nabla_\alpha\psi^\beta-\partial_\alpha\nabla_\gamma\psi^\beta)-(\partial_\alpha\dot{\Gamma}^\beta_\gamma-\dot{\Gamma}^\beta_\delta\Gamma^\delta_{\alpha\gamma})\psi^\gamma-(\dot{\Gamma}^\gamma_\alpha-\partial_\alpha V^\gamma)\nabla_\gamma\psi^\beta\\
\end{equation*}
\textrm{The last two expressions in parentheses can easily be simplified:}
\begin{equation*}
\psi^\gamma\partial_t\Gamma^\beta_{\alpha\gamma}-\Gamma^\beta_{\alpha\gamma}(\dot{\nabla}\psi^\gamma-\partial_t\psi^\gamma)-V^\gamma(\nabla_\gamma\nabla_\alpha\psi^\beta-\partial_\alpha\nabla_\gamma\psi^\beta)-(\nabla_\alpha\dot{\Gamma}^\beta_\gamma-\Gamma^\beta_{\alpha\delta}\dot{\Gamma}^\delta_\gamma)\psi^\gamma-(\Gamma^\gamma_{\alpha\delta}V^\delta-\tilde{C}B^\gamma_\alpha)\nabla_\gamma\psi^\beta\\
\end{equation*}
\textrm{The second term in parentheses can be simplified by recalling that the Riemann Curvature Tensor will appear if the indices of the two Covariant Derivatives can be reversed in the sense }$\nabla_\gamma\nabla_\alpha\psi^\beta=\nabla_\alpha\nabla_\gamma\psi^\beta-R^\beta_{\delta\alpha\gamma}\psi^\delta$\textrm{. Recalling the following identity }$\nabla_\alpha\nabla_\gamma\psi^\beta-\partial_\alpha\nabla_\gamma\psi^\beta=\Gamma^\beta_{\alpha\delta}\nabla_\gamma\psi^\delta-\Gamma^\delta_{\alpha\gamma}\nabla_\delta\psi^\beta$\textrm{ and collecting terms, the expression becomes reduced to:}
\begin{equation*}
\psi^\gamma\partial_t\Gamma^\beta_{\alpha\gamma}-\Gamma^\beta_{\alpha\gamma}(\dot{\nabla}\psi^\gamma-\partial_t\psi^\gamma)+V^\gamma R^\beta_{\phantom{.}\delta\alpha\gamma}\psi^\delta-V^\gamma\Gamma^\beta_{\alpha\delta}\nabla_\gamma\psi^\delta-(\nabla_\alpha\dot{\Gamma}^\beta_\gamma-\Gamma^\beta_{\alpha\delta}\dot{\Gamma}^\delta_\gamma)\psi^\gamma+\tilde{C}B^\gamma_\alpha\nabla_\gamma\psi^\beta
\end{equation*}
\textrm{Expanding the first expression in parentheses the expression is simplified to:}
\begin{equation*}
\psi^\gamma\partial_t\Gamma^\beta_{\alpha\gamma}-\Gamma^\beta_{\alpha\gamma}(\dot{\Gamma}^\gamma_\delta\psi^\delta-V^\delta\nabla_\delta\psi^\gamma)+V^\gamma R^\beta_{\phantom{.}\delta\alpha\gamma}\psi^\delta-V^\gamma\Gamma^\beta_{\alpha\delta}\nabla_\gamma\psi^\delta-(\nabla_\alpha\dot{\Gamma}^\beta_\gamma-\Gamma^\beta_{\alpha\delta}\dot{\Gamma}^\delta_\gamma)\psi^\gamma+\tilde{C}B^\gamma_\alpha\nabla_\gamma\psi^\beta
\end{equation*}
\textrm{Noticing that the first term in the first expression in parentheses annihilates the second term in the second expression in parentheses, and also noticing that the remaining }$V^\gamma$\textrm{ terms annihilate eachother, the entire Commutation Equation can be written as:}
\begin{equation*}
(\dot{\nabla}\nabla_\alpha-\nabla_\alpha\dot{\nabla})\psi^\beta=\psi^\gamma\partial_t\Gamma^\beta_{\alpha\gamma}+V^\gamma R^\beta_{\phantom{.}\delta\alpha\gamma}\psi^\delta-\psi^\gamma\nabla_\alpha\dot{\Gamma}^\beta_\gamma+\tilde{C}B^\gamma_\alpha\nabla_\gamma\psi^\beta
\end{equation*}
\textrm{Renaming the indices and grouping, the following is obtained:}
\begin{equation*}
(\dot{\nabla}\nabla_\alpha-\nabla_\alpha\dot{\nabla})\psi^\beta=(\partial_t\Gamma^\beta_{\alpha\gamma}+R^\beta_{\phantom{.}\gamma\alpha\delta}V^\delta-\nabla_\alpha\dot{\Gamma}^\beta_\gamma)\psi^\gamma+\tilde{C}B^\gamma_\alpha\nabla_\gamma\psi^\beta
\end{equation*}
\textrm{The opportunity will be taken to define a new differential object, the }\textbf{Temporal Curvature Tensor}\textrm{ \cite{BetterCMS} :}
\begin{equation}\label{TemCurvTen}
\dot{R}^\beta_{\phantom{.}\alpha\gamma}=\partial_t\Gamma^\beta_{\alpha\gamma}+R^\beta_{\phantom{.}\gamma\alpha\delta}V^\delta-\nabla_\alpha\dot{\Gamma}^\beta_\gamma
\end{equation}
\textrm{Which is so defined such that:}
\begin{equation}
(\dot{\nabla}\nabla_\alpha-\nabla_\alpha\dot{\nabla}-\tilde{C}B^\gamma_\alpha\nabla_\gamma)\psi^\beta=\dot{R}^\beta_{\phantom{.}\alpha\gamma}\psi^\gamma
\end{equation}
\textrm{The above relation can also be also be represented in commutation form as:}
\begin{equation}
([\dot{\nabla},\nabla_\alpha]-\tilde{C}B^\gamma_\alpha\nabla_\gamma)\psi^\beta=\dot{R}^\beta_{\phantom{.}\alpha\gamma}\psi^\gamma
\end{equation}
\textrm{It is worthwhile to note that }\textbf{no}\textrm{ term of }$\dot{R}^\beta_{\phantom{.}\alpha\gamma}$\textrm{ is a Tensor on its own (ie. }$\Gamma^\beta_{\alpha\gamma}$\textrm{, }$V^\gamma$\textrm{, and }$\dot{\Gamma}^\beta_\alpha$\textrm{), yet the entire expression must be a tensor by virtue of the commutating tensorial operators. This will be explored further in a later section.}
\subsection{The Temporal Connection on a Manifold}
\textrm{The above equation compared to Eq.(\ref{SpaTemInv}) mirror many of the properties that the original Riemann Curvature connection \cite{CMS, DifGeo1} shared in Eq.(\ref{SurfCommRiem}) compared to Eq.(\ref{SurfComm}).}\\\\\phantom{idnt}
\textrm{The properties can be more easily seen by observing the original definition of the Riemann Curvature connection. Considering two Vector Fields defined on a manifold, the connection is defined as:}
\begin{equation*}
\nabla_X\nabla_Y-\nabla_Y\nabla_X-\nabla_{\nabla_XY}
\end{equation*}
\textrm{Also expressable in commutation form as:}
\begin{equation*}
[\nabla_X,\nabla_Y]-\nabla_{\nabla_XY}
\end{equation*}
\textrm{On a manifold, allowing the two Vector Fields to be the components of the tangent space, the Vector Fields are orthogonal, and so the last term disappears \cite{DifGeo1}. In the Tensorial Time Derivative case }$\dot{\nabla}$\textrm{, the vector fields that Temporal Differentiation creates, }$\dot{\nabla}\textbf{T}$\textrm{ and tangent spatial vector fields, }$\textbf{X}=X^\alpha\textbf{S}_\alpha$\textrm{ are not necessarily orthogonal }\textbf{in the spatial sense}\textrm{, hence why the third term appears in the Tensorial Time Derivative Commutation:}
\begin{equation*}
[\dot{\nabla},\nabla_\alpha]-\tilde{C}B^\gamma_\alpha\nabla_\gamma
\end{equation*}
\textrm{Thus, it can be summarized that the commutation of the two derivative operators is:}
\begin{equation}
(\dot{\nabla}\nabla_\alpha-\nabla_\alpha\dot{\nabla})\psi^\beta=\dot{R}^\beta_{\phantom{.}\alpha\gamma}\psi^\gamma+\tilde{C}B^\gamma_\alpha\nabla_\gamma\psi^\beta
\end{equation}
\textrm{Similarly to the Riemann Curvature Tensor, using similar logic, it can be proven that for a covariant tensor defined on the surface:}
\begin{equation*}
([\dot{\nabla},\nabla_\alpha]-\tilde{C}B^\gamma_\alpha\nabla_\gamma)\psi_\beta=-\dot{R}^\gamma_{\phantom{.}\alpha\beta}\psi_\gamma
\end{equation*}
\textrm{As would expected for such a Curvature Tensor. This anti-symmetry is also observed for the regular Riemann Curvature Tensor \cite{CMS, TenAn2, DifGeo1}}
\subsection{Impossibility of `Temporally Curved' Static Surfaces}
\textrm{The Standard Riemann Curvature Tensor and Temporal Curvature Tensor differ strongly on one condition: For the Riemann Curvature Tensor, }\textbf{static surfaces }\underline{\textbf{have}}\textbf{ Riemann Curvature}\textrm{; this condition fails for the Temporal Curvature Tensor.}\\\\\phantom{idnt}
\textrm{This is of interest as the Tensor incorporates the Riemann Curvature Tensor, but in such a fashion that depends on time. This can be seen by the Temporal Curvature Tensor's Definition:}
\begin{equation*}
\dot{R}^\beta_{\phantom{.}\alpha\gamma}=\partial_t\Gamma^\beta_{\alpha\gamma}+R^\beta_{\phantom{.}\gamma\alpha\delta}V^\delta-\nabla_\alpha\dot{\Gamma}^\beta_\gamma
\end{equation*}
\textrm{If the Surface is }\textbf{static}\textrm{, then it can be intuitively stated that }$\textbf{V}=\textbf{0}$\textrm{; therefore this implies geometrically and can be proven easily analytically:}
\begin{equation*}
V^\beta=0\text{ and }\tilde{C}=0
\end{equation*}
\textrm{Thus, the last two terms vanish and the Temporal Curvature Tensor reduces to:}
\begin{equation*}
\dot{R}^\beta_{\phantom{.}\alpha\gamma}=\partial_t\Gamma^\beta_{\alpha\gamma}
\end{equation*}
\textrm{If the surface is static, then the Christoffel Symbols will depend entirely on the surface coordinates (ie. }$\Gamma^\beta_{\alpha\gamma}=\Gamma^\beta_{\alpha\gamma}(S)$\textrm{) and thus, its time partial derivative will vanish }$\partial_t\Gamma^\beta_{\alpha\beta}=0$\textrm{, producing the expected result for }\textbf{static}\textrm{ surfaces:}
\begin{equation*}
\dot{R}^\beta_{\phantom{.}\alpha\gamma}=0
\end{equation*}
\textrm{Thus this suggests that the Temporal Curvature Tensor not only is an appropiate indicator of Curvature, having incorporated the Riemann Curvature Tensor \cite{TenAn2}, but is also an appropiate indicator of the Dynamic Motion of a surface (having been defined using Surface Objects, and vanishing for }$\textbf{static}$\textrm{ surfaces).}
\subsection{Special Case: The Expanding Sphere \& Cylinder}
\textrm{It may appear from the equation that there should be Temporal Curvature for every non ricci-flat smoothly deforming surface, but this is not the case as shown with two prototypical examples below:}
\subsubsection*{The Expanding Sphere}
\textrm{It can be shown that for the Sphere with a radius of }$R(t)$\textrm{, its Temporal Curvature can be simplified greatly. This can be demonstrated by analyzing the Sphere's differential objects on the surface. For the sphere, it has a vanishing Surface Speed, }$V^\alpha=\textbf{V}\cdot\textbf{S}_\alpha$\textrm{. This makes sense geometrically and can be proven analytically:}
\begin{align*}
&\phantom{equation}\text{\underline{Since the Sphere is always expanding in the direction of its normal}}\\
&\phantom{equation}\text{\underline{thus, the projection of its velocity onto the tangent space should vanish.}}
\end{align*}
\textrm{Therefore in this case, the Temporal Curvature Tensor reduces to:}
\begin{equation*}
\dot{R}^\beta_{\phantom{.}\alpha\gamma}=\partial_t\Gamma^\beta_{\alpha\gamma}+\nabla_\alpha(\tilde{C}B^\beta_\gamma)
\end{equation*}
\textrm{For a Sphere (which is at all points equivalent), the Christoffel Symbols are not dependent on the Radius }$R(t)$\textrm{, and therfore are not dependent on time. In addition for the Sphere, neither the Surface Normal Velocity }$C=R'(t)$\textrm{ nor  the Mean Curvature Tensor, }$B^\beta_\gamma=-R^{-1}\delta^\beta_\gamma$\textrm{ are dependent on position }$(\theta, \phi)$\textrm{. So it can be determined that for the Sphere, }$\partial_t\Gamma^\beta_{\alpha\gamma}=0$\textrm{ and for the Covariant Derivative term, utilizing the properties of the Covariant Derivative and the properties of the Sphere:}
\begin{equation*}
\nabla_\alpha(\tilde{C}B^\beta_\gamma)=\nabla_\alpha\left(-\frac{R'(t)}{R(t)}\delta^\beta_\gamma\right)=-\frac{R'(t)}{R(t)}\nabla_\alpha\delta^\beta_\gamma=0
\end{equation*}
\textrm{This establishes that for the Sphere, the Temporal Curvature Tensor }$\textbf{vanishes}$\textrm{!}
\begin{equation}
\dot{R}^\beta_{\phantom{.}\alpha\gamma}=0
\end{equation}
\subsubsection*{The Expanding Cylinder}
\textrm{Examining a Cylinder with a radius of }$R(t)$\textrm{, its Temporal Curvature can be simplified greatly. This can be demonstrated by analyzing the Cylinder's differential objects on the surface. For the Cylinder, it has a vanishing Surface Speed, }$V^\alpha=\textbf{V}\cdot\textbf{S}_\alpha$\textrm{. This makes sense geometrically and can be proven analytically:}
\begin{align*}
&\phantom{equation}\textrm{Since the Cylinder is always expanding in the direction of its normal}\\
&\phantom{equation}\textrm{thus, the projection of its velocity onto the tangent space should vanish.}
\end{align*}
\textrm{Therefore, in the Cylinder's case, the Temporal Curvature Tensor reduces to:}
\begin{equation*}
\dot{R}^\beta_{\phantom{.}\alpha\gamma}=\partial_t\Gamma^\beta_{\alpha\gamma}+\nabla_\alpha(\tilde{C}B^\beta_\gamma)
\end{equation*}
\textrm{For a Cylinder, the Christoffel Symbols are not only independent of time, they also vansh for the cylinder! }$\Gamma^\gamma_{\alpha\beta}=0$\textrm{. In addition for the Cylinder, neither the Surface Normal Velocity }$C=R'(t)$\textrm{ nor  the Mean Curvature Tensor, (}$B^\beta_\gamma=\text{Diag(}-1/R\text{, 0)}$\textrm{) are dependent on position }$(\theta, z)$\textrm{. So it can be determined that for the Cylinder, }$\partial_t\Gamma^\beta_{\alpha\gamma}=0$\textrm{ and for the Covariant Derivative term, since the Christoffel Symbols vanish, the result is immediately obtained:}
\begin{equation}
\dot{R}^\gamma_{\alpha\beta}=0
\end{equation}
\textrm{Thus, for the expanding Cylinder \& Sphere, the Temporal Curvature Tensor }$\textbf{vanishes}$\textrm{!}
\subsection{Simplifying the Temporal Curvature Tensor}
\textrm{Noticing for the expanding Sphere \& Cylinder that the Temporal Curvature Tensor entirely vanishes, it can be inquired if this result can observed by simplifying the Tensor further. So, motivated by the previous two examples the Tensor will be attempted to be simplified. First, the first term of the Tensor, }$\partial_t\Gamma^\beta_{\alpha\gamma}$\textrm{ is decomposed into the Christoffel Symbol of the First Kind }$\Gamma_{\alpha\gamma\beta}$\textrm{, and the Metric Tensor }$S_{\alpha\beta}$\textrm{, then the product rule is applied the product rule:}
\begin{equation*}
\partial_t\Gamma^\beta_{\alpha\gamma}=\Gamma_{\alpha\gamma\delta}(\partial_tS^{\beta\delta})+S^{\beta\delta}(\partial_t\Gamma_{\alpha\gamma\delta})
\end{equation*}
\textrm{Then, recalling the formula }$\Gamma_{\alpha\gamma\delta}=1/2(\partial_\gamma S_{\alpha\delta}+\partial_\alpha S_{\gamma\delta}-\partial_\delta S_{\alpha\gamma})$\textrm{, the expansion of the Christoffel Symbol of the First Kind is performed and simplified:}
\begin{equation*}
\partial_t\Gamma^\beta_{\alpha\gamma}=\Gamma_{\alpha\gamma\delta}(\partial_t\partial S^{\beta\delta})+\frac{1}{2}S^{\beta\delta}\partial_t\left(\partial_\gamma S_{\alpha\delta}+\partial_\alpha S_{\gamma\delta}-\partial_\delta S_{\alpha\gamma}\right)
\end{equation*}
\textrm{Knowing that the partial time derivative }$\partial_t$\textrm{ and surface partial derivatives }$\partial_\alpha$\textrm{ commute, the time derivative is expanded and replaced with the partial spatial derivatives:}
\begin{equation*}
\partial_t\Gamma^\beta_{\alpha\gamma}=\Gamma_{\alpha\gamma\delta}(\partial_t S^{\beta\delta})+\frac{1}{2}S^{\beta\delta}\left(\partial_\gamma\partial_tS_{\alpha\delta}+\partial_\alpha\partial_tS_{\gamma\delta}-\partial_\delta\partial_tS_{\alpha\gamma}\right)
\end{equation*}
\textrm{It is at this point that the Metric Time Derivative, }$K_{\alpha\beta}$\textrm{ is defined as:}
\begin{equation}
K_{\alpha\beta}=\partial_tS_{\alpha\beta}
\end{equation}
\textrm{This can be expanded in terms of the Normal and Surface Speed, by Rearranging the Tensorial Time Derivative definition of the Metric Tensor, and remembering that }$\dot{\nabla}S_{\alpha\beta}=0${ and that }$\nabla_\gamma S_{\alpha\beta}=0$\textrm{:}
\begin{equation*}
\dot{\nabla}S_{\alpha\beta}=\partial_tS_{\alpha\beta}-V^\gamma\nabla_\gamma S_{\alpha\beta}-\dot{\Gamma}^\gamma_{\alpha}S_{\gamma\alpha}-\dot{\Gamma}^\gamma_{\beta}S_{\gamma\alpha}\phantom{.}\to\phantom{.}\partial-tS_{\alpha\beta}=K_{\alpha\beta}=\dot{\Gamma}^\gamma_{\alpha}S_{\gamma\beta}+\dot{\Gamma}^\gamma_{\beta}S_{\gamma\alpha}
\end{equation*}
\textrm{Using the Symmbetry of the Mean Curvature Tensor }$B_{\alpha\beta}$\textrm{, the Metric's Time Derivative can be expressed as another form:}
\begin{equation}
K_{\alpha\beta}=\nabla_\alpha V_\beta+\nabla_\beta V_\alpha-2\tilde{C}B_{\alpha\beta}
\end{equation}
\textrm{By the identical method, one can obtain as expected:}
\begin{equation}
\frac{\partial S^{\alpha\beta}}{\partial t}=-K^{\alpha\beta}
\end{equation}
\textrm{Using these definitions, the first term of the Tensor, }$\partial_t\Gamma^\beta_{\alpha\gamma}$\textrm{ can be simplified further:}
\begin{equation*}
\partial_t\Gamma^\beta_{\alpha\gamma}=-K^{\beta\delta}\Gamma_{\alpha\gamma\delta}+\frac{1}{2}S^{\beta\delta}\left(\partial_\gamma K_{\alpha\delta}+\partial_\alpha K_{\gamma\delta}-\partial_\delta K_{\alpha\gamma}\right)
\end{equation*}
\textrm{For each partial derivative of the Metric Time Derivative, }$\partial_\alpha K_{\beta\gamma}$\textrm{ it will be expressed as a Pseudo-Covariant Derivative }$\nabla_\alpha K_{\beta\gamma}$\textrm{. Since each of the Metric Time Derivatives contain two indices, by the Covariant Derivative, there will be }$\textbf{six}$\textrm{ Christoffel Symbols of the Second Kind which appear by the rearrangement.}\\\\\phantom{idnt}
\textrm{Allowing the symmetry of the bottom two indices of the Christoffel Symbols }$\Gamma^\gamma_{[\alpha \beta]}=0$\textrm{, the two Christoffel Symbols generated by }$\partial_\delta K_{\alpha\gamma}$\textrm{, cancel out two of the Christoffel Symbols in the first two bracket terms (ie. }$\partial_\gamma K_{\alpha\delta}\text{ \& }\partial_\alpha K_{\gamma\delta}$\textrm{). Thus only two Christoffel Symbols remain and three Covariant Derivatives:}
\begin{equation*}
\partial_t\Gamma^\beta_{\alpha\gamma}=-K^{\beta\delta}\Gamma_{\alpha\gamma\delta}+\frac{1}{2}S^{\beta\delta}\left(\nabla_\gamma K_{\alpha\delta}+\nabla_\alpha K_{\gamma\delta}-\nabla_\delta K_{\alpha\gamma}+\Gamma^\epsilon_{\alpha\gamma}K_{\epsilon\delta}+\Gamma^\epsilon_{\gamma\alpha}K_{\epsilon\delta}\right)
\end{equation*}
\textrm{Once more, utilizing the symmetry of the Christoffel Symbols, the following is obtained:}
\begin{equation*}
\partial_t\Gamma^\beta_{\alpha\gamma}=-K^{\beta\delta}\Gamma_{\alpha\gamma\delta}+\frac{1}{2}S^{\beta\delta}\left(\nabla_\gamma K_{\alpha\delta}+\nabla_\alpha K_{\gamma\delta}-\nabla_\delta K_{\alpha\gamma}+2\Gamma^\epsilon_{\alpha\gamma}K_{\epsilon\delta}\right)
\end{equation*}
\textrm{Two implications are worth noting. If the action of the Contravariant Metric Tensor is allowed on the }$2\Gamma^\epsilon_{\alpha\gamma}K_{\epsilon\delta}$\textrm{ term in brackets raising the final index producing }$2\Gamma^\epsilon_{\alpha\gamma}K^\beta_\epsilon$\textrm{ while lowering the index on the Christoffel Symbols making it of the First Kind }$2\Gamma_{\alpha\gamma\epsilon}K^{\epsilon\beta}$\textrm{, and expanding the 1/2 factor, the first term and last term in the expression annihilate. This results in the Simplification:} 
\begin{equation*}
\partial_t\Gamma^\beta_{\alpha\gamma}=\frac{1}{2}S^{\beta\delta}\left(\nabla_\gamma K_{\alpha\delta}+\nabla_\alpha K_{\gamma\delta}-\nabla_\delta K_{\alpha\gamma}\right)
\end{equation*}
\textrm{There are similarities between this definition, and the definition of the Christoffel Symbol of the First Kind. This equation can be simplified greatly recalling that }$K_{\alpha\beta}=\nabla_\alpha V_\beta+\nabla_\beta V_\alpha-2\tilde{C}B_{\alpha\beta}$\textrm{. After expanding all the Metric Time Derivatives according to this definition there remain several Twice-Iterated Pseudo-Covariant Derivatives of the Surface Speed }$\nabla_\alpha\nabla_\beta V_\gamma$\textrm{. Many of these can be simplified into the definition of the Riemann Curvature Tensor. After expanding and simplifying the expression using the definition of the Riemann Curvature tensor, the following is obtained:}
\begin{equation*}
\frac{\partial}{\partial t}\Gamma^\beta_{\alpha\gamma}=\frac{1}{2}S^{\beta\delta}\big(2\nabla_\alpha\nabla_\gamma V_\delta-V_\epsilon(R^\epsilon_{\phantom{.}\delta\gamma\alpha}+R^\epsilon_{\phantom{.}\gamma\alpha\delta}+R^\epsilon_{\phantom{.}\alpha\gamma\delta})+2\nabla_\delta(\tilde{C}B_{\alpha\gamma})-2\nabla_\alpha(\tilde{C}B_{\gamma\delta})-2\nabla_\gamma(\tilde{C}B_{\alpha\delta})\big)
\end{equation*}
\textrm{The expression in brackets can be simplified greatly. The Bianchi Identity states a cyclic relation among the Riemann Curvature Tensors \cite{CMS}:}
\begin{equation*}
R^\epsilon_{\phantom{.}\gamma\alpha\delta}+R^\epsilon_{\phantom{.}\delta\gamma\alpha}+R^\epsilon_{\phantom{.}\alpha\delta\gamma}=0\phantom{..}\to\phantom{..}R^\epsilon_{\phantom{.}\gamma\alpha\delta}+R^\epsilon_{\phantom{.}\delta\gamma\alpha}=-R^\epsilon_{\phantom{.}\alpha\delta\gamma}
\end{equation*}
\textrm{Recalling that reversing the last two indices of the Riemann Curvature Tensor makes the Tensor negative \cite{CMS}:}
\begin{equation*}
R^\epsilon_{\phantom{.}\gamma\alpha\delta}+R^\epsilon_{\phantom{.}\delta\gamma\alpha}=R^\epsilon_{\phantom{.}\alpha\gamma\delta}
\end{equation*}
\textrm{Thus, the Time Differentiation of the Christoffel Symbol of the Second Kind is simplified to:}
\begin{equation*}
\partial_t\Gamma^\beta_{\alpha\gamma}=\frac{1}{2}S^{\beta\delta}\big(2\nabla_\alpha\nabla_\gamma V_\delta-2V_\epsilon R^\epsilon_{\phantom{.}\alpha\gamma\delta}+2\nabla_\delta(\tilde{C}B_{\alpha\gamma})-2\nabla_\alpha(\tilde{C}B_{\gamma\delta})-2\nabla_\gamma(\tilde{C}B_{\alpha\delta})\big)
\end{equation*}
\textrm{The Riemann Curvature's index is then lowered to its covariant form which raises the Surface Speed's index and thus, the 1/2 factor may finally be expanded into the brackets to obtain:}
\begin{equation*}
\partial_t\Gamma^\beta_{\alpha\gamma}=S^{\beta\delta}\Big(\nabla_\alpha\nabla_\gamma V_\delta-V^\epsilon R_{\epsilon\alpha\gamma\delta}+\nabla_\delta(\tilde{C}B_{\alpha\gamma})-\nabla_\alpha(\tilde{C}B_{\gamma\delta})-\nabla_\gamma(\tilde{C}B_{\alpha\delta})\Big)
\end{equation*}
\textrm{After this, the Riemann Curvature Tensor's Indices are permuted, and the Metric Tensor is expanded into the brackets to obtain:}
\begin{equation*}
\partial_t\Gamma^\beta_{\alpha\gamma}=\nabla_\alpha\nabla_\gamma V^\beta-V^\epsilon R^\beta_{\phantom{.}\gamma\alpha\epsilon}+\nabla^\beta(\tilde{C}B_{\alpha\gamma})-\nabla_\alpha(\tilde{C}B^\beta_\gamma)-\nabla_\gamma(\tilde{C}B^\beta_\alpha)
\end{equation*}
\textrm{The first and fourth term can be grouped together to form:}
\begin{equation*}
\partial_t\Gamma^\beta_{\alpha\gamma}=\nabla_\alpha\dot{\Gamma}^\beta_\gamma-V^\epsilon R^\beta_{\phantom{.}\gamma\alpha\epsilon}+\nabla^\beta(\tilde{C}B_{\alpha\gamma})-\nabla_\gamma(\tilde{C}B^\beta_\alpha)
\end{equation*}
\textrm{This is the CMS definition of the Time Derivative of the Christoffel Symbols of the Second Kind. The beauty of the simplifcation appears upon review with the whole definition of the Temporal Curvature Tensor. Recalling the Temporal Curvature Tensor's definition from earlier in Eq.(\ref{TemCurvTen}), \& substituting this Christoffel Symbol's derivative into the Temporal Curvature Tensor; the first and second term dissapear, and the Temporal Curvature Tensor reduces to the simplification:}
\begin{equation}
\dot{R}^\beta_{\alpha\gamma}=\nabla^\beta(\tilde{C}B_{\alpha\gamma})-\nabla_\gamma(\tilde{C}B^\beta_\alpha)
\end{equation}
\textrm{In this fashion, the need for the Christoffel Symbols and Riemann Curvature Tensor is eliminated. The Temporal Curvature is now dependent on essentially }$\textbf{V}$\textrm{ (from which }$\tilde{C}$\textrm{ is obtained) and }$\textbf{S}_\alpha$\textrm{ (from which }$\textbf{N}$\textrm{ and therefore }$B_{\alpha\beta}$\textrm{ is obtained). The Temporal Curvature Tensor may actually be simplified further. The product is expanded using the product rule, the following develops:}
\begin{equation*}
\dot{R}^\beta_{\alpha\gamma}=\tilde{C}\nabla^\beta B_{\alpha\gamma}-\tilde{C}\nabla_\gamma B^\beta_\alpha+B_{\alpha\gamma}\nabla^\beta\tilde{C}-B^\beta_\alpha\nabla_\gamma\tilde{C}
\end{equation*}
\textrm{Lowering first two terms using the metric tensor, the following term develops in parentheses:}
\begin{equation*}
\dot{R}^\beta_{\alpha\gamma}=\tilde{C}S^{\beta\delta}(\nabla_\delta B_{\gamma\alpha}-\nabla_\gamma B_{\delta\alpha})+B_{\alpha\gamma}\nabla^\beta\tilde{C}-B^\beta_\alpha\nabla_\gamma\tilde{C}
\end{equation*}
\textrm{Recalling the Codazzi Equation \cite{CMS}, the two terms in the parentheses annihilate eachother, resulting in:}
\begin{equation}\label{FinalTemCurvTen}
\dot{R}^\beta_{\alpha\gamma}=B_{\alpha\gamma}\nabla^\beta\tilde{C}-B^\beta_\alpha\nabla_\gamma\tilde{C}
\end{equation}
\textrm{This is a final reduced expression of the Temporal Curvature Tensor and is significant in that comparing it to the earlier definition of the Tensor, }\textbf{every term is in fact Tensorial}\textrm{.}
\subsection{The Action of the Spatio-Temporal Commutation on several Differential Objects}
\textrm{Though Eq.(\ref{FinalTemCurvTen}) has a long derivation to obtain, it can also appear in a much quicker fashion. Analyzing the commutation of the operators }$\dot{\nabla}\text{ and }\nabla_\alpha$\textrm{ on the tangent vectors }$\textbf{S}_\alpha$\textrm{, the commutation takes the following form:}
\begin{equation*}
(\dot{\nabla}\nabla_\alpha-\nabla_\alpha\dot{\nabla})\textbf{S}_\beta=-\dot{R}^\gamma_{\alpha\beta}\textbf{S}_\gamma+\tilde{C}B^\gamma_\alpha\nabla_\gamma\textbf{S}_\beta
\end{equation*}
\textrm{Using the definitions of the Tensorial Time Derivative and Covariant Derivatives on the tangent vectors (recalling }$\dot{\nabla}\textbf{S}_\alpha=\textbf{N}\nabla_\alpha C$\textrm{ and }$\nabla_\alpha\textbf{S}_\beta=\textbf{N}B_{\alpha\beta}$\textrm{), and factoring all terms between the Tangent Vectors and the Normal, the following develops:}
\begin{equation*}
\textbf{N}(\dot{\nabla}B_{\alpha\beta}-\nabla_\alpha\nabla_\beta\tilde{C}-\tilde{C}B^\gamma_\alpha B_{\gamma\beta})+\textbf{S}_\gamma(\dot{R}^\gamma_{\alpha\beta}-B_{\alpha\beta}\nabla^\gamma\tilde{C}+B^\gamma_\alpha\nabla_\beta\tilde{C})=0
\end{equation*}
\textrm{We already recognize the normal projection of the equation as the Invariant Time Derivative of the Curvature Tensor in Eq.(\ref{CMStable}); this is analogous to the Codazzi Equation \cite{CMS} in the sense that it is the Tensorial Time Derivative analogue of }$\textbf{N}\cdot(\nabla_\beta\nabla_\gamma-\nabla_\gamma\nabla_\beta)\textbf{S}_\alpha$\textrm{ relating the Tensorial Time Derivative of the Curvature Tensor to itself much like the Codazzi Equation related the derivatives of the Curvature Tensor.}\\\\\phantom{idnt}
\textrm{Though recently introduced, the Tangential Projection of the Equation can be identified with the recent introduction of the simplification of the Temporal Curvature Tensor; this relation is analogous to Gauss' Theorema Egregium in the sense that is the Tensorial Time Derivative analogue of }$\textbf{S}^\delta\cdot(\nabla_\beta\nabla_\gamma-\nabla_\gamma\nabla_\beta)\textbf{S}_\alpha$\textrm{ \cite{CMS}. It relates the Temporal Curvature Tensor to the Curvature Tensor and Normal Speed of the Surface much like the Theorema Egregium related the Riemann Curvature Tensor to the Curvature Tensor of the Surface. This confirms the simplification of the Temporal Curvature Tensor, and the Time Derivative of the Christoffel Symbols all in one fluid factorization. In analyzing the Normal in a similar fashion using the }$[\dot{\nabla},\nabla_\alpha]$\textrm{ commutation, the following develops:}
\begin{equation*}
\textbf{N}(B_{\alpha\gamma}\nabla^\gamma\tilde{C}-\nabla_\gamma\tilde{C}B^\gamma_\alpha)+\textbf{S}_\gamma(\nabla_\alpha\nabla^\gamma\tilde{C}+\tilde{C}B^\beta_\alpha B^\gamma_\beta-\dot{\nabla}B^\gamma_\alpha)=0
\end{equation*}
\textrm{In this case, the result is in-fact a statement of two statements which agree with previously established relations. The Tangential Projection of this Equation is the Equavilent of the Normal Projection from the previous Commutation of }$\textbf{S}_\alpha$\textrm{ which is already familiar; the Normal Projection of this Equation implies another concept entirely stated as:}
\begin{equation}\label{CommRel1}
B_{\alpha\gamma}\nabla^\gamma\tilde{C}-B^\gamma_\alpha\nabla_\gamma\tilde{C}=0
\end{equation}
\textrm{Is an observation which will be adressed later.}
\subsection{Differential of the Curvature and Temporal Tensor}
\textrm{It can be seen that the disappearance of the Temporal Curvature Tensor on the Sphere and Cylinder was due to the time-independence of their Christoffel Symbols and the fact that their Curvatures were not dependent on the position.}\\\\\phantom{idnt}
\textrm{One value can be obtained from this Tensor as its }$\textbf{trace}$\textrm{. Much like the Tensor itself vanished, so will its trace. The trace gives a counterintuitive notion that while the expanding surfaces were in fact curved }$\textbf{and}$\textrm{ dynamic, the trace of their Temporal Curvature Tensor vanished. Motivated by this example, the trace is writ in a tensorial form:}
\begin{equation}
\tilde{\delta}^\tau\hat{\textbf{R}}=\dot{R}^\beta_{\alpha\beta}\textbf{S}^\alpha
\end{equation}
\textrm{A constant may also be obtained from this value the }\textbf{Scalar Temporal Curvature Differential}\textrm{:}
\begin{equation}
\tilde{\delta}^\tau\kappa=|\tilde{\delta}^\tau\hat{\textbf{R}}|
\end{equation}
\textrm{Expanding the value in terms of components, the formula is expressed as:}
\begin{equation}
\tilde{\delta}^\tau\kappa=\sqrt{S^{\alpha\beta}\dot{R}^\gamma_{\alpha\gamma}\dot{R}^\sigma_{\beta\sigma}}
\end{equation}
\textrm{At present moment, the only statement which can be made with certainty is that all uniformly expanding surfaces (such as the Expanding Cylinder \& Sphere) observe the following law:}
\begin{equation*}
\tilde{\delta}^\tau\kappa=0
\end{equation*}
\textrm{Though the definition does not appears to be invalid, nor does it describe much, an observation can be made }\textbf{in general}\textrm{ about surfaces. Utilizing Eq.(\ref{FinalTemCurvTen}) and taking its trace, the following relation is obtained:}
\begin{equation*}
\dot{R}^\beta_{\alpha\beta}=B_{\alpha\beta}\nabla^\beta\tilde{C}-B^\beta_\alpha\nabla_\beta\tilde{C}
\end{equation*}
\textrm{Since the indices can be raised/lowered at whim, we notice that for }\underline{\textbf{all surfaces}}\textrm{, this Tensor Vanishes stating the following relation:}
\begin{equation}
\dot{R}^\beta_{\alpha\beta}=0
\end{equation}
\textrm{ This equation applied for all surfaces, not just for expanding cylinders and spheres. This equation is essentially a confirmation of the relation obtained through the Commutation in Eq.(\ref{CommRel1}). In addition, it can be stated that for }$\underline{all}$\textrm{ surfaces:}
\begin{equation}
\tilde{\delta}^\tau\kappa=0
\end{equation}
\subsection{Ambient Temporal Curvature Tensor}
\textrm{For the operator }$\left(\nabla_\alpha\nabla_\beta-\nabla_\beta\nabla_\alpha\right)$\textrm{ it has been shown that when applied to an surface tensor with }\underline{\textbf{ambient}}\textrm{ indices of order (1,0) or (0,1) as in }$\vec{\psi}=\psi^i\textbf{Z}_i=\psi_i\textbf{Z}^i$\textrm{, the operator }\textbf{does}\textrm{ commute (unlike the surface index case above):}
\begin{equation*}
\left(\nabla_\alpha\nabla_\beta-\nabla_\beta\nabla_\alpha\right)\psi^i=\left(\nabla_\alpha\nabla_\beta-\nabla_\beta\nabla_\alpha\right)\psi_i=0
\end{equation*}
\textrm{This effect is due to the Tensorial Time Derivative of the Ambient Base vanishing }$\dot{\nabla}\textbf{Z}_i$\textrm{, but the Surface Basis does not, }$\dot{\nabla}\textbf{S}_\alpha=\textbf{N}\nabla_\alpha\tilde{C}$\textrm{. It can be inquired if a similar effect occurs with the newly defined operator }$(\dot{\nabla}\nabla_\alpha-\nabla_\alpha\dot{\nabla})$\textrm{. The result after using all conventions is that:}
\begin{equation*}
(\dot{\nabla}\nabla_\alpha-\nabla_\alpha\dot{\nabla})\psi^i=\left[\partial_t\left(Z^j_{\phantom{.}\alpha}\Gamma^i_{jk}\right)-\nabla_{\alpha}\left(V^j\Gamma^i_{jk}\right)-V^\beta R^i_{\phantom{.}klj}Z^j_{\phantom{.}\alpha}Z^l_{\phantom{.}\beta}\right]\psi^k+\tilde{C}B^\beta_\alpha\nabla_\beta\psi^i
\end{equation*}
\textrm{For a covariant ambient tensor, it can also be obtained then:}
\begin{equation*}
(\dot{\nabla}\nabla_\alpha-\nabla_\alpha\dot{\nabla})\psi_i=-\left[\partial_t\left(Z^j_{\phantom{.}\alpha}\Gamma^k_{ij}\right)-\nabla_{\alpha}\left(V^j\Gamma^k_{ij}\right)-V^\beta R^k_{\phantom{.}ilj}Z^j_{\phantom{.}\alpha}Z^k_{\phantom{.}\beta}\right]\psi_k+\tilde{C}B^\beta_\alpha\nabla_\beta\psi_i
\end{equation*}
\textrm{Based on these equations, an Ambient Temporal Curvature Tensor may be defined as the following:}
\begin{equation}
\dot{R}^i_{\phantom{.}\alpha k}=\partial_t\left(Z^j_{\phantom{.}\alpha}\Gamma^k_{ij}\right)-\nabla_{\alpha}\left(V^j\Gamma^k_{ij}\right)-V^\beta R^i_{\phantom{.}klj}Z^j_{\phantom{.}\alpha}Z^l_{\phantom{.}\beta}
\end{equation}
\textrm{So defined such that:}
\begin{equation*}
(\dot{\nabla}\nabla_\alpha-\nabla_\alpha\dot{\nabla}-\tilde{C}B^\beta_\alpha\nabla_\beta)\psi^i=\dot{R}^i_{\phantom{.}\alpha k}\psi^k
\end{equation*}
\textrm{As expected from such commutations, switching the tensor to a covariant one, much like earlier, the familiar result appears:}
\begin{equation*}
(\dot{\nabla}\nabla_\alpha-\nabla_\alpha\dot{\nabla}-\tilde{C}B^\beta_\alpha\nabla_\beta)\psi_i=-\dot{R}^k_{\phantom{.}\alpha i}\psi_k
\end{equation*}
\textrm{Using the Ambient Temporal Curvature Tensor, the equations may be re-defined to read as:}
\begin{equation}
\begin{matrix}
(\dot{\nabla}\nabla_\alpha-\nabla_\alpha\dot{\nabla})\psi^i=\dot{R}^i_{\phantom{.}\alpha k}\psi^k+\tilde{C}B^\beta_\alpha\nabla_\beta\psi^i \\
(\dot{\nabla}\nabla_\alpha-\nabla_\alpha\dot{\nabla})\psi_i=-\dot{R}^k_{\phantom{.}\alpha i}\psi_k+\tilde{C}B^\beta_\alpha\nabla_\beta\psi_i
\end{matrix}
\end{equation}
\textrm{Interestingly, for a Euclidean Space, }$(x,y,z)$\textrm{, it can be easily found that }$\Gamma^i_{jk}=0$\textrm{. As such, the conclusion is reached that for a Euclidean Space:}
\begin{equation}
\dot{R}^i_{\phantom{.}\alpha k}=0
\end{equation}
\textrm{After work, the Ambient Temporal Curvature Tensor can be abbreviated into a tensorial form as the following:}
\begin{equation}
\dot{R}^i_{\alpha k}=R^i_{\phantom{.}klj}Z^l_{\phantom{.}\alpha}CN^j
\end{equation}
\textrm{This form highlights the dependency of the Commutation on the Curvature of the Space. In this case, for whichever surface, if the surface is static, the Ambient Temporal Curvature Tensor }$\textbf{will}$\textrm{ vanish.}
\section{The Tri-Velocity Equation}
\textrm{Though the table given in Eq.(\ref{CMStable}) is well diverse, there are several relations which are omitted. Three objects which are omitted from the table are the Ambient Speed, }$\textbf{V}$\textrm{, Surface Speed, }$V^\alpha$\textrm{ and the Surface Velocity, }$\tilde{C}$\textrm{ All three velocities are related in the following statement \cite{CMS}:}
\begin{equation}
\textbf{V}=\tilde{C}\textbf{N}+\textbf{S}_\alpha V^\alpha
\end{equation}
\textrm{The Equation states that the Ambient Speed, }$\textbf{V}$\textrm{ (which is the speed of a surface in ambient space) is equivalent to its Normal Speed, given by }$\tilde{C}\textbf{N}$\textrm{, and its tangential Surface Speed, }$\textbf{S}_\alpha V^\alpha$\textrm{ combined. The relation is interesting in that since the Normal Field and Shift Tensor obey the relations:}
\begin{equation*}
\left\{\begin{matrix}
\textbf{N}\cdot\textbf{N}=1 \\
\textbf{S}_\alpha\cdot\textbf{S}^\beta=\delta^{\beta}_{\phantom{.}\alpha} \\
\textbf{N}\cdot\textbf{S}_\alpha=0
\end{matrix}\right.
\end{equation*}
\textrm{Dotting the above relation with the dual Surface Basis (}$\textbf{S}^\beta$\textrm{), or with the Normal Field (}$\textbf{N}$\textrm{) yield the appropiate relations, respectively:}
\begin{equation*}
\left\{V^\beta=\textbf{V}\cdot\textbf{S}^\beta\phantom{.},\phantom{.}\tilde{C}=\textbf{V}\cdot\textbf{N}\right\}
\end{equation*}
\textrm{Though not a Tensor as per CMS \cite{CMS}, the Ambient Speed }$\textbf{V}$\textrm{ is a Pseudotensor in that it posesses sufficient structure to observe Differentiation; thus, Tensorial Operators such as the Tensorial Time Derivative }$\dot{\nabla}$\textrm{ may still be performed on }$\textbf{V}$\textrm{ so long as it is acknowledged that the result }$\textbf{may}$\textrm{ not be a tensor. Thus, the Ambient Speed's Pseudotensorial Time Derivative may be stated as:}
\begin{align*}
\dot{\nabla}\textbf{V}=\dot{\nabla}(\tilde{C}\textbf{N})+\dot{\nabla}(\textbf{S}_{\alpha}V^\alpha)
\end{align*}
\textrm{Using the relations from the table in Eq.(\ref{CMStable}), the following is the derivative's decomposition:}
\begin{align*}
\dot{\nabla}\textbf{V}=\textbf{N}\dot{\nabla}\tilde{C}-\textbf{S}_{\alpha}\tilde{C}\nabla^\alpha\tilde{C}+\textbf{S}_\alpha\dot{\nabla}V^\alpha+V^\alpha \textbf{N}\nabla_\alpha\tilde{C}
\end{align*}
\textrm{Factoring this into its normal and tangential projections, the }\textbf{Tri-Velocity Equation}\textrm{ may be obtained:}
\begin{equation}\label{TriVEq}
\dot{\nabla}\textbf{V}=\textbf{N}(\dot{\nabla}\tilde{C}+V^\alpha\nabla_\alpha\tilde{C})+\textbf{S}_{\alpha}(\dot{\nabla}V^\alpha-\tilde{C}\nabla^\alpha\tilde{C})
\end{equation}
\textrm{The Tri-Velocity Equation is so named as Eq.(\ref{TriVEq}) relates the Pseudotensorial Time Derivatives of all three measures of speed: }$\textbf{V}\text{, }\tilde{C}\text{, and }V^\alpha$\textrm{. If the Ambient Speed's }$\textbf{V}$\textrm{ Temsorial Time Derivative is known, then the surface velocity's }$\tilde{C}$\textrm{ can be determined and the surface speed's }$V^\gamma$\textrm{ can also be known by recalling that }$\tilde{C}=\textbf{V}\cdot\textbf{N}$\textrm{ and }$V^\alpha=\textbf{S}^\alpha\cdot\textbf{V}$\textrm{. Applying the Tensorial Time Derivative Operator to these relations, the following two equations may be obtained:}
\begin{equation}
\dot{\nabla}\tilde{C}=\textbf{N}\cdot\dot{\nabla}\textbf{V}-V^\alpha\nabla_\alpha\tilde{C}\phantom{...},\phantom{...}\dot{\nabla}V^\alpha=\textbf{S}^\beta\cdot\dot{\nabla}\textbf{V}+\tilde{C}\nabla^\alpha\tilde{C}
\end{equation}
\textrm{Thus, all three measures of speed are linked togather by Eq.(\ref{TriVEq}).}
\subsection{Temporal Acceleration of Surface}
\textrm{The Acceleration of the Surface is a Non-Tensorial Vector quantity defined by:}
\begin{equation}
\textbf{A}=\partial_t\textbf{V}
\end{equation}
\textrm{Where }$\textbf{V}=V^i\textbf{Z}_i$\textrm{. Thus, the equation may be expanded in component-form:}
\begin{equation*}
\textbf{A}=(\partial_tV^i)\textbf{Z}_i+V^i(\partial_t\textbf{Z}_i)
\end{equation*}
\textrm{Applying the chain rule and allowing the introduction of the Christoffel Symbols, the following expression develops:}
\begin{equation}
\partial_t\textbf{Z}_i=(\partial_j\textbf{Z}_i)(\partial_tZ^j)=\Gamma^k_{ij}V^j\textbf{Z}_k
\end{equation}
\textrm{Using the above simplification, the second term may be simplified, and thus, the whole equation to the following:}
\begin{equation*}
\textbf{A}=\left(\partial_tV^k+V^iV^j\Gamma^k_{ij}\right)\textbf{Z}_k
\end{equation*}
\textrm{If we define the acceleration component by }$A^i$\textrm{, such that }$\textbf{A}=A^i\textbf{Z}_i$\textrm{, then we see:}
\begin{equation*}
A^k=\partial_tV^k+V^iV^j\Gamma^k_{ij}
\end{equation*}
\subsection{The Ambient Tensor Acceleration}
\textrm{By applying the definition of the Invariant Time Derivative to the Ambient Velocity, we obtain:}
\begin{equation*}
\dot{\nabla}\textbf{V}=\partial_t\textbf{V}-V^\alpha\nabla_\alpha\textbf{V}
\end{equation*}
\textrm{This can be abbreviated into the acceleration tensor:}
\begin{equation*}
\dot{\nabla}\textbf{V}=\textbf{A}-V^\alpha\nabla_\alpha\textbf{V}
\end{equation*}
\textrm{We expand the covariant derivative by definition to the following:}
\begin{equation*}
\dot{\nabla}\textbf{V}=\textbf{A}-V^\alpha\partial_\alpha\textbf{V}
\end{equation*}
\textrm{This term can be simplified by considering the first term in brackets. We see that:}
\begin{equation*}
\partial_\alpha\textbf{V}=\partial_\alpha\partial_t\textbf{R}=\partial_t\partial_\alpha\textbf{R}=\partial_t\textbf{S}_\alpha
\end{equation*}
\textrm{We also notice that we can substitute the equivalency:}
\begin{equation*}
\partial_t\textbf{S}_\alpha=\dot{\nabla}\textbf{S}_\alpha+V^\beta\nabla_\beta\textbf{S}_\alpha+\dot{\Gamma}^\beta_\alpha\textbf{S}_\beta=(\nabla_\alpha\tilde{C}+V^\beta B_{\alpha\beta})\textbf{N}+\dot{\Gamma}^\beta_\alpha\textbf{S}_\beta
\end{equation*}
\textrm{Thus, we see:}
\begin{equation*}
\dot{\nabla}\textbf{V}=\textbf{A}-V^\alpha(\nabla_\alpha\tilde{C}+V^\beta B_{\alpha\beta})\textbf{N}-V^\alpha\dot{\Gamma}^\beta_\alpha\textbf{S}_\beta
\end{equation*}
\textrm{We can use the Tri-Velocity equation to obtain:}
\begin{equation}
\dot{\nabla}\tilde{C}=\textbf{N}\cdot\textbf{A}-2V^\alpha\nabla_\alpha\tilde{C}-V^\alpha V^\beta B_{\alpha\beta}
\end{equation}
\textrm{And also:}
\begin{equation}
\dot{\nabla}V^\gamma=\textbf{S}^\gamma\cdot\textbf{A}-V^\alpha\dot{\Gamma}^\gamma_\alpha+\tilde{C}\nabla^\gamma\tilde{C}
\end{equation}
\textrm{These are essential to have. Even though the Tri-Veloctiy Equation consitutes a useful relation between the Ambient Velocity, Surface Speed, and Surface Velocity, the above equations explicitly state the relation for each unique Speed.}
\subsection{Three Speeds, Three Tensorial Accelerations}
\textrm{Since the operator is an Invariant time derivative for Invariants, and a proper Tensorial Time Derivative for Tensors \cite{CMS}, the above equations also constitute a sort of Tensorial Accelerations, the Ambient Acceleration Tensor:}
\begin{equation}
\mathbb{A}^i=A^i-N^i(V^\alpha\nabla_\alpha\tilde{C}+V^\alpha V^\beta B_{\alpha\beta})-V^\alpha\dot{\Gamma}^\beta_\alpha Z^i_{\phantom{.}\beta}
\end{equation}
\textrm{, the Surface Acceleration Tensor:}
\begin{equation}
\tilde{\mathcal{A}}^\gamma=\dot{\nabla}V^\gamma=Z^\gamma_{\phantom{.}i}A^i-V^\alpha\dot{\Gamma}^\gamma_\alpha+\tilde{C}\nabla^\gamma\tilde{C}
\end{equation}
\textrm{, and the Invariant Normal Acceleration:}
\begin{equation}
\hat{\mathcal{A}}=\dot{\nabla}C=N_iA^i-2V^\alpha\nabla_\alpha\tilde{C}-V^\alpha V^\beta B_{\alpha\beta}
\end{equation}
\textrm{Also, by this, we notice that the Tri-Velocity Equation can be expressed as a `Tri-Acceleration' Equation:}
\begin{equation}
\mathbb{A}^i=N^i(\hat{\mathcal{A}}+V^\alpha\nabla_\alpha\tilde{C})+Z^i_{\phantom{.}\alpha}(\tilde{\mathcal{A}}^\alpha-\tilde{C}\nabla^\alpha\tilde{C})
\end{equation}
\textrm{Interestingly, this demonstrates that the Normal Projection of the Invariant Normal Acceleration and the Tangential Projection of the Surface Acceleration Tensor is not enough to represent the acceleration of an arbitrary surface in ambient space:}
\begin{equation*}
\mathbb{A}^i-N^i\hat{\mathcal{A}}-Z^i_{\phantom{.}\alpha}\tilde{\mathcal{A}}^\alpha=N^iV^\alpha\nabla_\alpha\tilde{C}-Z^i_{\phantom{.}\alpha}\tilde{C}\nabla^\alpha\tilde{C}
\end{equation*}
\textrm{If we define the following commutation:}
\begin{equation}
(\hat{\mathcal{N}}_*\tilde{\circ}\mathcal{T}_*)^{ij}_\alpha=N^iZ^j_{\alpha}-Z^i_{\phantom{.}\alpha}N^j
\end{equation}
\textrm{And define the component of the Ambient Acceleration not represented in the Surface/Normal Projections of the Surface/Normal Accelerations as }$\tilde{\Delta}\hat{\mathcal{A}}^i$\textrm{, Then, we can abbreviate the Tri-Acceleration Equation as:}
\begin{equation}
\tilde{\Delta}\hat{\mathcal{A}}^i=(\hat{\mathcal{N}}_*\tilde{\circ}\mathcal{T}_*)^{ij}_\alpha V_j\nabla^\alpha\tilde{C}
\end{equation}
\textrm{Finally the Commutator can be seen to have interesting properties:}
\begin{equation}
\left\{\begin{matrix}
&\frac{1}{2}\epsilon_{ijk}\epsilon^{\alpha\beta}(\hat{\mathcal{N}}_*\tilde{\circ}\mathcal{T}_*)^{ij}_\alpha=Z^\beta_k \\
&Z_{ij}(\hat{\mathcal{N}}_*\tilde{\circ}\mathcal{T}_*)^{ij}_\alpha=0
\end{matrix}\right.
\end{equation}
\subsection{Practical Applications}
\subsubsection*{Kinetic Energy}
\textrm{In many cases, we see that the Kinetic Energy of a Surface is given by \cite{CMS}:}
\begin{equation*}
\mathcal{K}=\frac{1}{2}\rho_S\int_S|\textbf{V}|^2dS
\end{equation*}
\textrm{This equation can be expressed tensorially as:}
\begin{equation}
\mathcal{K}=\frac{1}{2}\rho_S\int_SV_iV^idS
\end{equation}
\textrm{We recognize that the rate of kinetic energy change is given by:}
\begin{equation*}
\frac{d\mathcal{K}}{dt}=\frac{1}{2}\rho_S\frac{d}{dt}\int_SV_iV^idS
\end{equation*}
\textrm{We can also simplify this using the known derivative of a surface integral:}
\begin{equation*}
\frac{d\mathcal{K}}{dt}=\frac{1}{2}\rho_S\int_S(\dot{\nabla}-\tilde{C}B^\alpha_\alpha)V_iV^idS
\end{equation*}
\textrm{Therefore, we can simplify the equation as such:}
\begin{equation*}
\frac{d\mathcal{K}}{dt}=\frac{1}{2}\rho_S\int_S\dot{\nabla}(V_iV^i)-\tilde{C}B^\alpha_\alpha V_iV^idS
\end{equation*}
\textrm{and can simplify the first term in the Integral using our definition of the Ambient Acceleration:}
\begin{equation*}
\frac{d\mathcal{K}}{dt}=\frac{1}{2}\rho_S\int_S2V_i\mathbb{A}^i-\tilde{C}B^\alpha_\alpha V_iV^idS
\end{equation*}
\textrm{Expanding the product in the last term, the Equation simplifies to:}
\begin{equation*}
\frac{d\mathcal{K}}{dt}=\frac{1}{2}\rho_S\int_S2V_i\mathbb{A}^i-\tilde{C}B^\alpha_\alpha(\tilde{C}^2+V^\beta V_\beta)dS
\end{equation*}
\textrm{And finally expand the 1/2 factor into the whole integrand to obtain:}
\begin{equation}
\frac{d\mathcal{K}}{dt}=\rho_S\int_SV_i\mathbb{A}^i-\frac{1}{2}\tilde{C}B^\alpha_\alpha(\tilde{C}^2+V^\beta V_\beta)dS
\end{equation}
\subsubsection*{Power}
\textrm{We also know that the power of a moving surface is given by:}
\begin{equation}
\mathcal{P}=\rho_S\int_S\vec{\alpha}\cdot\textbf{V}dS
\end{equation}
\textrm{Where }$\vec{\alpha}$\textrm{ is the acceleration felt by a point mass in a corresponding field, such that }$\vec{\alpha}=\frac{1}{m}\textbf{F}$\textrm{. This integral form can also be put into a tensorial form:}
\begin{equation*}
\mathcal{P}=\rho_S\int_SV_i\alpha^idS
\end{equation*}
\subsubsection*{Work-Energy Theorem}
\textrm{Surfaces and Particles alike abide by the principle of Work-Energy. This states that in a surface in which all energy given by External Forces and there is null net working, the surface will convert all the energy into kinetic energy. This theorem is well established and has undergone several extensions over the past 50 years \cite{WorkEn}. This is stated usually as:}
\begin{equation}
\mathcal{P}=\frac{d\mathcal{K}}{dt}
\end{equation}
\textrm{Based on the following principle, both sides of the equation were derived above. In its integral form, we have:}
\begin{equation*}
\rho_S\int_SV_i\alpha^idS=\rho_S\int_SV_i\mathbb{A}^i-\frac{1}{2}\tilde{C}B^\alpha_\alpha(\tilde{C}^2+V^\beta V_\beta)dS
\end{equation*}
\textrm{If we simplify this and assume that an vanishing integral over an arbitrary domain implies a vanishing integrand, then we find the following equation of motion:}
\begin{equation}
V_i\alpha^i=V_i\mathbb{A}^i-\frac{1}{2}\tilde{C}B^\alpha_\alpha(\tilde{C}^2+V^\beta V_\beta)
\end{equation}
\textrm{This equation can also be rewrote as:}
\begin{equation}
V_i(\mathbb{A}^i-\alpha^i)=\frac{1}{2}\tilde{C}B^\alpha_\alpha(\tilde{C}^2+V^\beta V_\beta)
\end{equation}
\textrm{Or in its vector form:}
\begin{equation}
\textbf{V}\cdot(\vec{\mathbb{A}}-\vec{\alpha})=\frac{1}{2}\tilde{C}B^\alpha_\alpha(\tilde{C}^2+V^\beta V_\beta)
\end{equation}
\textrm{This equation has several observations. In the case that the term on the right dissapears, then the equation implies that:}
\begin{equation*}
\vec{\mathbb{A}}=\vec{\alpha}
\end{equation*}
\textrm{This special case suggests that the acceleration of each point on the surface will experience the same acceleration as a point mass would in a field. We have seen in the special cases of expanding surface that they have no tangential velocity, this reducing the equations to the following form:}
\begin{equation*}
N_i(A^i-\alpha^i)=\frac{1}{2}C^2B^\alpha_\alpha
\end{equation*}
\textrm{If the surface is moving slowly, we assume that }$\tilde{C}^2\approx0$\textrm{ and the equations are reduced to}
\begin{equation*}
\mathbb{A}^i=\alpha^i
\end{equation*}
\section{Conclusion}
\textrm{The Calculus of Moving Surfaces, CMS,  is such a rich budding field of Calculus which both the Mathematics discipline and Physics discipline can embrace. It offers interesting insights into the next intuitive step from Differential Geometry, extending it to Dynamics Manifolds and also offers exciting new opportunities and models for Theoretical Physics as it pertains to model surface phenomena, an area of research which is of increasing importance. \cite{Clo2DSurf, EMagSurf} }\\\\\phantom{idnt}
\textrm{While it presents a few preliminary short-comings in imediately admitting similarities with classical Physics, further work in deriving analogues of theorems within Classical Mechanics in a CMS framework will prove to be fruitful. CMS provides interesting insights in models ranging from Soap-Films \cite{CMS, Micelles}, to Biological Membranes \cite{Micelles}, to classical optimization problems encountered in elementary Multivariable Calculus \cite{Droplets}, deriving theorems and relations such as the Young-Laplace Equation, Hadamard Principle, Minimal Surface Problem, with ease \cite{CMS, BetterCMS}. In addition, the sub-discipline also finds utility in solving problems of static surfaces by assuming that the time-dependence allows a dynamics manifold to evolve to its equilibrium as an optimized static manifold based on an energy integral \cite{BetterCMS}.}\\\\\phantom{idnt}
\textrm{Its young nature only implies that the more investigation that is done in this beautiful sub-discipline of Differential Geometry, the deeper insights that will be discovered into the Nature of Surfaces and the more opportunities for modelling which will be uncovered.}

\end{document}